\documentclass[12pt]{article}
\usepackage{amssymb,amsmath,amsfonts,amssymb}
\newcommand{\be}{\begin{equation}}
\newcommand{\ee}{\end{equation}}
\newcommand{\la}{\label}
\newcommand{\ba}{\begin{array}{c}}
\newcommand{\ea}{\end{array}}
\newcommand{\Rm}{{\mathbb R}}

\newcommand{\nax}{\nabla_x}

\newcommand{\pak}{\partial_k}
\newcommand{\pai}{\partial_i}
\newcommand{\paj}{\partial_j}

\newcommand{\sij}{\sigma_{ij}}
\newcommand{\sip}{\sigma_{ip}}
\newcommand{\skp}{\sigma_{kp}}
\newcommand{\sjp}{\sigma_{jp}}

\newcommand{\dx}{{\text{div}}_x}
\newcommand{\dd}{{\mathcal{D}}}

\newcommand{\qv}[2]{\left\langle #1, #2 \right\rangle}
\newcommand{\at}[1]{\big|_{#1}}
\newcommand{\sub}[1]{^{(#1)}}

\newtheorem{thm}{Theorem}
\newtheorem{prop}{Proposition}
\newtheorem{lemma}{Lemma}

\title{Stochastic Lagrangian Transport\\ and\\ Generalized Relative Entropies}
\author{Peter Constantin\\ Department of Mathematics, 
The University of Chicago\\ 5734 S. University Avenue, Chicago, Illinois 60637\
\and
Gautam Iyer\\Department of Mathematics, Stanford University\\
Bldg.\ 380, 450 Serra Mall, Stanford, CA 94304
}

\begin{document}
\maketitle

\begin{abstract} {We discuss stochastic representations of advection diffusion equations with variable diffusivity, stochastic integrals of motion and generalized relative entropies.}
\end{abstract}


\noindent{\bf AMS - MSC numbers: 35K45, 60H30.}

\section{Introduction}

Recently, Michel, Mischler and Perthame \cite{Perthame1} discovered a remarkable property of certain unstable linear equations, in which decay of relative entropies takes place. Their observation was applied to population dynamics models, but the list of applications is growing. Of course, relative entropies have been used for a long time in kinetic theory and conservation laws. However, the decay of relative entropies, was known before 
only in stable, self-adjoint situations in which a global attracting steady solution exists and no flow advection is present \cite{villani}. The property of decay of relative entropies was slightly  generalized to variable diffusion coefficients and applied to Smoluchowski systems in \cite{c}. A stochastic interpretation and proof in the case of constant diffusion coefficients was given in \cite{c1}. Here we provide a stochastic interpretation and proof in
the case of variable diffusion coefficients. The method of proof and
concepts are of more general interest \cite{ci,gautam}.

We consider a linear operator
\be
{\dd}\rho = \nu\pai(a_{ij}\paj \rho) - \dx(U\rho) + V\rho\la{opdiv}
\ee
in ${\mathbb R}^n$, where
\be
U(x,t) = \left (U_j(x,t)\right )_{j=1,\dots n}\la{U}
\ee
is a smooth ($C^2$) function, $V = V(x,t)$ is a continuous and bounded 
scalar potential and
\be
a_{ij}(x,t) = \sip(x,t)\sjp(x,t)\la{aij}
\ee
with the matrix
\be
\sigma(x,t) = (\sij(x,t))_{ij}
\la{sigma}
\ee
a given smooth ($C^2$) matrix. We assume that $\sigma$ is bounded and $U$ and $\nax \sigma$ decay at infinity.
We use the shorthand notation $A(D)$ for the operator
\be 
A(D)\rho = a_{ij}\pai\paj\rho
\la{ad}
\ee
and use also the non-divergence form
\be
{\dd}\rho = \nu A(D)\rho - u\cdot\nax \rho + P\rho
\la{op}
\ee
where
\be
u_j(x,t) = U_j(x,t) - \nu\pai(a_{ij}(x,t))\la{u}
\ee
and
\be
P = V - \dx(U).
\la{P}
\ee
The formal adjoint of the operator $\dd$  
in $L^2({\mathbb R}^n)$ is
\be
{\dd}^*\phi = \nu\pai(a_{ij}\paj\phi) + U \cdot\nax \phi + V\phi \la{ddad}.
\ee
The following is the result of Michel, Mischler and Perthame:
\begin{thm}\cite{Perthame,Perthame1}
Let $f$ be a solution of
\be
\partial_t f = {\dd}f
\la{feq}
\ee
and let $\rho >0$ be a positive solution of the same equation,
\be
\partial_t \rho = {\dd}\rho. 
\la{rhoeq}
\ee
Let $H$ be a smooth convex function of one variable and let
$\phi$ be a non-negative function obeying pointwise
\be
\partial_t \phi + {\dd}^*\phi = 0 \la{phieq}.
\ee
Then
\be
\frac{d}{dt}\int H\left(\frac{f}{\rho}\right) \phi\rho dx \le 0.
\la{in}
\ee
\end{thm}

\section{Stochastic Lagrangian Flow}
In order to represent solutions of equations like (\ref{feq}) we consider
the drift
\be
v_j(x,t) = u_j + 2\nu (\pak\sjp)\skp = U_j -\nu (\pak\skp)\sjp + \nu(\pak\sjp)\skp.\la{vj}
\ee
Let $X(a,t)$ be the strong solution of the stochastic differential system
\be
dX_j(t) = v_j(X,t) dt + \sqrt{2\nu}\sjp(X,t) dW_p\la{dx}
\ee
with initial data
\be
X(a,0) = a.\la{aid}
\ee
Here $W$ is a standard Brownian process in ${\mathbb R}^n$ starting at time zero from the origin. This process will be fixed throughout the paper and all
measurability issues will be with respect to the filtration associated to it and all almost sure statements will be with respect to the probability measure
on the standard Wiener space.
We will need the following result:

\begin{thm}\la{Athm}
The inverse of the flow map $a\mapsto X(a,t)$, the stochastic map
\be
x\mapsto A(x,t)\la{a}
\ee
exists almost surely and satisfies its defining relations
$$
X(A(x,t),t) = x, \,\,\ \forall x\in {\mathbb R}^n, \quad A(X(a,t),t) = a, \,\,\, \forall a\in{\mathbb R}^n, \,\,\, \forall t, \, a.s.
$$
The map $X$ is smooth and the determinant
\be
D(a,t) = \det{(\partial_a X(a,t))}\la{d}
\ee
obeys the SDE
\be
\ba
d(\det{(\partial_a X(a,t))} = [\det(\partial_a X(a,t))]\times\\\left\{\left [ (\dx v)(x,t)+ 2\nu E(x,t)\right ]_{|x= X(a,t)}dt + \sqrt{2\nu} (\pak(\skp))(x,t)_{|x= X(a,t)}dW_p\right\}
\ea
\la{deteq}
\ee
with
\be
E(x,t) = \sum_{i<j}\sum_p \det(\pai\sjp)_{ij}.
\la{e}
\ee
The map $A(x,t)$ satisfies the stochastic partial differential system
\be
dA_j + \left (u\cdot\nax A_j - \nu A(D) A_j\right)dt + \sqrt{2\nu}(\pak A_j)\skp dW_p = 0\la{aeq}
\ee 
with initial data
$$
A(x,0) = 0.
$$
\end{thm}

\noindent{\bf Remark.} In the statement above, $\det(\pai\sjp)_{ij}$
refers to the determinant of the two-by-two matrix $(\partial_r\skp)$ with $r,k\in \{i, j\}$ for fixed $i<j$ and $p$. Theorem \ref{Athm} was originally proved in \cite{ci} for constant coefficients and in \cite{gautam} for variable coefficients. For completeness, we reproduce the proof (with variable coefficients, as stated above) in Appendix \ref{AthmPf}.

\section{Stochastically Passive Scalars and \\Feynman-Kac Formula }

We consider first deterministic smooth time-independent functions $f_0$ and note that the
functions $\theta = \theta_{f_0} (x,t) = f_0(A(x,t))$ are stochastically passive in the sense that they obey the equation
\be
d\theta + \left( u\cdot\nax \theta - \nu A(D)\theta\right)dt + \sqrt{2\nu}\pak\theta\skp dW_p = 0\la{thetaeq}
\ee
with initial data
\be
\theta(x,0) = f(x).\la{thetaid}
\ee
Solutions of the SPDE (\ref{thetaeq}) form an algebra; in particular, products 
of solutions are solutions, a nontrivial fact due to the presence of the stochastic term. The expected values of these scalars obey advection-diffusion equations and do not form an algebra in general, if $\nu>0$.
We consider now the function
\be
I(a,t) = \exp{\left\{\int_0^t P(X(a,s),s)ds\right\}}
\la{i}
\ee
where $P(x,t)$ is given in (\ref{P}) and consider the function
\be
\psi =\psi_{f_0}(x,t) = \theta_{f_0}(x,t) I(A(x,t),t)\la{psif}
\ee
We have
\begin{thm}
The process $\psi = \psi_{f_0}$ given by 
\be
\psi(x,t) = f_0(A(x,t))\exp{\left \{\int_0^t P(X(a,s),s)ds_{|a= A(x,t)}\right\}}\la{psi}
\ee
solves
\be
d\psi - \left (\dd\psi \right)dt + {\sqrt{2\nu}}\nax\psi \sigma  dW = 0
\la{psieq}
\ee
with initial datum $\psi(x,0) =f_0(x)$.
\end{thm}
The proof of this result follows using stochastic calculus \cite{Ka}, \cite{Ku}. Indeed, the function $I(a,t)$
obeys
\be
\partial_t I(a,t) = P(X(a,t),t) I(a,t)
\la{ieq}
\ee
pathwise (almost surely).  Then, a calculation using ({\ref{aeq}}) (see \cite{ci}, \cite{gautam})  
shows that the function
\be
J(x,t) = I(A(x,t),t)
\la{J}
\ee
solves 
\be
dJ + (u\cdot\nax J - PJ - \nu A(D)J)dt + \sqrt{2\nu}\nax J\sigma dW = 0.
\la{jeq}
\ee
The function $\psi_{f_0}$ is the product
$$
\psi_{f_0} = \theta_{f_0} J,
$$
and therefore, from It\^{o}'s formula
$$
d\psi_{f_0} = Jd\theta + \theta dJ + d\langle J,\theta\rangle
$$
and the equations obeyed by $J$, $\theta$, we have
$$
\ba
d\psi_{f_0} = \\
 (-u\cdot\nax\psi_{f_0} + P\psi_{f_0} + \nu JA(D)\theta + \nu\theta A(D) J
 + 2\nu (\pak J)\skp (\paj \theta)\sjp )dt \\ -\sqrt{2\nu}\nax \psi_{f_0}\sigma dW.
\ea
$$
This means
$$
d\psi = (-u\cdot\nax \psi + P\psi +\nu A(D)\psi)dt -\sqrt{2\nu}\nax \psi_{f_0}\sigma\
 dW.
$$
Because of (\ref{op}) we have (\ref{psieq}).

\section{Stochastic Integrals of Motion.}
\begin{prop}
Consider a deterministic function $\phi$ that solves (\ref{phieq}). Then the function
\be
M(a,t) = \phi(X(a,t),t)\det{(\partial_a X(a,t))}\exp\left\{\int_0^t P(X(a,s),s)ds\right\}
\la{M}
\ee
is a martingale. 
\end{prop}
\noindent{\bf Proof.}
We start by writing 
$$
M(a,t) = \Phi (a,t)I(a,t) D(a,t)
$$
with
$$
\Phi(a,t) = \phi(X(a,t),t),
$$
$I$ given above in (\ref{i}) and $D$ given in (\ref{d}). 
Next, we compute the equation obeyed by $\Phi I$. In view of (\ref{ieq}) and using It\^{o}'s formula we have
$$
\ba
d(\Phi I) = I\left \{\left (\partial_t\phi(X(a,t),t) + P\phi(X(a,t),t) \right)dt \right. +\\ + \left.  \nax\phi_{|X(a,t)}\cdot dX + \frac{1}{2}\pai\paj \phi_{|X(a,t)} d\langle X_i, X_j\rangle\right\},
\ea
$$
which gives, in view of (\ref{dx})
$$
\ba
d(\Phi I) = \\ =
\left\{ \partial_t\phi + P\phi + v\cdot\nax \phi + \nu A(D)\phi\right\}_{|X(a,t)}dt + \sqrt{2\nu}I((\pai\phi)\sip)_{|X(a,t)}dW_p.
\ea
$$
Using (\ref{phieq}) and (\ref{vj}) we have
\be
\ba
d(\Phi I) =  I\left\{-2\nu(\pak(\skp))\sjp (\paj\phi) - (\dx U)\phi\right\}_{|X(a,t)}dt + \\ +\sqrt{2\nu}I((\paj \phi)\sjp)_{|X(a,t)}dW_p.
\ea
\la{phiieq}
\ee
Now, by It\^{o},
$$
dM = Dd(\Phi I) + \Phi I D + d\langle D, \Phi I\rangle.
$$
In view of (\ref{deteq}) and (\ref{phiieq}) we have
$$
d\langle D, \Phi I\rangle = 2\nu D I\left\{(\pak\skp)\sjp (\paj\phi)\right\}_{|X(a,t)}dt
$$
and consequently the terms $\pm2\nu DI(\pak\skp)\sjp (\paj\phi) dt$ cancel and we obtain
$$
\ba
dM = ID\left\{- (\dx U)\phi + (\dx v + 2\nu E)\right\}_{|X(a,t)}dt +\\ + \sqrt{2\nu}ID \left \{\sjp (\paj\phi) + (\pak\skp)\right \}_{|X(a,t)}dW_p.
\ea
$$
Now, in view of (\ref{vj}) we have that
$$
(\dx v) - (\dx U) = \nu \paj[(\pak\sjp)\skp]- \nu \paj[(\pak\skp)\sjp]
$$ 
and therefore the coefficient of $dt$ in $dM$ is
$$
DI\left \{2\nu E + \nu \paj[(\pak\sjp)\skp]- \nu \paj[(\pak\skp)\sjp]\right\}
$$ 
Now
$$
\ba
DI\left\{\nu \paj[(\pak\sjp)\skp]- \nu \paj[(\pak\skp)\sjp]\right\} =\\
DI\left\{\nu(\pak\sjp)(\paj\skp) -\nu(\pak\skp)(\paj\sjp)\right\}=\\=
DI2\sum_{{k<j}}\sum_p\left\{\nu(\pak\sjp)(\paj\skp)- \nu(\pak\skp)(\paj\sjp)\right\}=\\ = -2\nu E
\ea
$$
and therefore the coefficient of $dt$ in $dM$ vanishes. We obtained
\be
dM =  \sqrt{2\nu}ID \left \{\sjp (\paj\phi) + (\pak\skp)\right \}_{|X(a,t)}dW_p\la{Meq},
\ee
that is, $M$ is the martingale
$$
\ba
M(a,t) = \\
\phi(a,0) + \sqrt{2\nu}\int_0^t I(a,s)D(a,s)\left \{\sjp (\paj\phi) + (\pak\skp)\right \}_{|X(a,s)}dW_p(s)
\ea
$$

\begin{thm} Let $h_0$ and $\rho_0$ be smooth time independent deterministic functions. Consider the stochastically passive scalar $\theta_{h_0}(x,t) = h_0(A(x,t))$ and the process $\psi_{\rho_0}$ of (\ref{psi}) with initial datum  $\rho_0$. Consider also $\phi (x,t)$, a deterministic solution of (\ref{phieq}). Then the random variable
\be
{\mathcal {E}}(t) = \int_{{\mathbb R}^n}\phi(x,t)\psi_{\rho_0}(x,t)\theta_{h_0}(x,t)dx
\la{sten}
\ee
is a martingale.
In particular
\be
{\mathbb E}({\mathcal {E}}(t)) = \int_{{\mathbb R}^n}\phi(a,0)\rho_0(a)h_0(a)da
\la{conse}
\ee
holds.
\end{thm}
\noindent{\bf Proof.} In view of the change of variables formula and the definition of $\psi_{\rho_0}$ we have that
\be
{\mathcal{E}}(t) = \int_{{\mathbb R}^n} M(a,t) \rho_0(a)h_0(a)da
\la{repe}
\ee
with $M$ given in (\ref{M}). The result follows then from the previous proposition. More precisely
\be
\ba
d{\mathcal {E}} = \\
\sqrt{2\nu}\left\{\int_{{\mathbb{R}^n}}\exp{\left\{\int_0^s P(X(a,\tau),\tau)d\tau_{|a= A(x,s)}\right\}} \left \{\sjp (\paj\phi) + (\pak\skp)\right \}dx\right\}dW_p
\ea
\la{Eeq}
\ee
gives explicitly the SDE obeyed by ${\mathcal {E}}$.

\section{Generalized Relative Entropies}
We take now a smooth deterministic, time independent function $H$ of one variable, a deterministic solution of (\ref{phieq}), two smooth deterministic, time independent functions $f_0$ and $\rho_0$, of which $\rho_0$ is strictly positive. We form the processes $\psi_{\rho_0}$ and $\psi_{f_0}$ given by the expressions (\ref{psi}). Then it t follows that
$$
\ba
\psi_{\rho_0}(x,t)\phi(x,t)H\left(\frac{\psi_{f_0}(x,t)}{\psi_{\rho_{0}}(x,t)}\right)
= \\ 
\psi_{\rho_0}(x,t)\phi(x,t)H\left(\frac{f_0(A(x,t))}{\rho_0(A(x,t))}\right )
\ea
$$
holds. Thus, the quantity of interest, $\psi_{\rho_0}\phi H\left(\frac{\psi_{f_0}}{\psi_{\rho_0}}\right)$,  is the product of a stochastically passive scalar, $\psi_{\rho_0}$ and $\phi$. By the previous theorem  we have that
\be
{\mathcal{E}}(t) = \int\psi_{\rho_{0}}(x,t)H\left (\frac{\psi_{f_0}(x,t)}{\psi_{\rho_{0}}(x,t)}\right )\phi(x,t) dx
\la{int}
\ee
is a martingale. The expected value is then constant in time:
\be
\frac{d}{dt}{\mathbb {E}}\left\{\int \psi_{\rho_0}H\left (\frac{\psi_{f_0}}{\psi_{\rho_0}}\right)\phi dx\right\} = 0.
\la{cons}
\ee
If we denote
\be
f(x,t) = {\mathbb{E}}\psi_{f_0}(x,t)
\la{f}
\ee
and
\be
\rho(x,t) = {\mathbb {E}}\psi_{\rho_0}(x,t)
\la{rho}
\ee
we have from (\ref{psieq}) that $f$ solves (\ref{feq}), $\rho>0$ solves
 (\ref{rhoeq}). We prove that we have (\ref{in}). 

The starting point is (\ref{cons}). In view of (\ref{f}) and (\ref{rho}), the
statement that needs to be proved is
\be
\int {\mathbb{E}}\left(\psi_{\rho_0}\right) H\left (\frac{{\mathbb{E}}(\psi_{f_0})}{{\mathbb{E}}(\psi_{\rho_0})}\right )\phi dx\le {\mathbb {E}}\left\{\int 
\psi_{\rho_0} H\left (\frac{\psi_{f_0}}{\psi_{\rho_0}}\right)\phi dx \right \}
\la{to}
\ee 
The conservation (\ref{cons}) works for any $H$, but we expect (\ref{to})
to hold only for convex $H$. Indeed, (\ref{to}) can be reduced to a Jensen
inequality. We claim more, that for all $x,t$ we have
\be
{\mathbb{E}}\left(\psi_{\rho_0}\right) H\left (\frac{{\mathbb{E}}(\psi_{f_0})}{{\mathbb{E}}(\psi_{\rho_0})}\right ) \le  {\mathbb {E}}\left\{
\psi_{\rho_0} H\left (\frac{\psi_{f_0}}{\psi_{\rho_0}}\right)\right\}
\la{tod}
\ee
Considering the functions 
\be
g = \frac{\psi_{\rho_0}}{{\mathbb{E}}(\psi_{\rho_0})}\la{g}
\ee
and
\be
v= \frac{\psi_{f_0}}{{\mathbb {E}}(\psi_{\rho_0})}\la{v}
\ee
we see that (\ref{tod}) becomes
\be
H\left ({\mathbb{E}}(v)\right) \le {\mathbb {E}}\left\{g H\left (\frac{v}{g}\right)\right\}.\la{todo}
\ee
This, however, is nothing but Jensen's inequality for the probability measure
$$
Ph = {\mathbb{E}}(gh),
$$
$$
H\left (P\left(\frac{v}{g}\right)\right ) \le PH\left (\frac{v}{g}\right ).
$$

\appendix
\section{Proof of Theorem \ref{Athm}}\label{AthmPf}
We devote this appendix to proving Theorem \ref{Athm}. The original proof can be found in \cite{ci} for constant coefficients, and in \cite{gautam} for variable coefficients.

\begin{lemma}\la{detAlemma}
Let $X$ be the stochastic flow defined by \eqref{dx}, \eqref{aid}. Then the map $X$ is spatially smooth (almost surely), and the determinant $D = \det( \nabla X)$ satisfies the equation
$$
dD = D \left[ \left(\nabla \cdot v + 2 \nu E \right)\,dt + \sqrt{2\nu} \pak\skp \,dW_p \right]
$$
where
$$
E = \tfrac{1}{2}\left[ \partial_i \sigma_{ip} \partial_j \sigma_{jp} - \partial_j\sigma_{ip} \partial_i\sigma_{jp} \right] .
$$
\end{lemma}
\noindent{\bf Proof.} Differentiating (\ref{dx}) we have
\be
d(\partial_a X_j) = \pak v_j \partial_a X_k\,dt + \sqrt{2\nu}\pak\sjp \partial_a X_k\,dW_p.\la{dparx}
\ee
Let $S^n$ be the permutation group on $n$ symbols, and $\epsilon_\tau$ denote the signature of the permutation $\tau \in S^n$. By It\^o's formula,
\begin{align}
dD  &= \sum_{\substack{\tau \in S^n\\b=1\dots n}} \epsilon_\tau \bigg[ \prod_{c\neq b} \partial_c X_{\tau_c} \,d\left( \partial_b X_{\tau_b} \right) + \sum_{c<b} \prod_{d \neq b,c} \partial_d X_{\tau_d} \,d\qv{\partial_b X_{\tau_b}}{\partial_c X_{\tau_c}}\bigg]	\nonumber\\
    &= \sum_{\substack{\tau \in S^n\\b=1\dots n}} \epsilon_\tau \bigg[ \partial_k v_{\tau_b} \partial_b X_k \prod_{c\neq b} \partial_c X_{\tau_c} \,dt + \sqrt{2\nu} \partial_k \sigma_{\tau_b,p} \partial_b X_k \prod_{c\neq b} \partial_c X_{\tau_c} \, dW_p + \nonumber\\
    &\qquad\qquad\qquad + \nu\sum_{c \neq b} \partial_b X_l \partial_l \sigma_{\tau_b,p} \, \partial_c X_m \partial_m \sigma_{\tau_c,p} \prod_{d\neq b,c} \partial_d X_{\tau_d} \, dt \bigg] \la{D1}
\end{align}
We compute each of the terms above individually:
\begin{align*}
\sum_{\substack{\tau \in S^n\\b=1\dots n}} \epsilon_\tau \partial_k v_{\tau_b} \partial_b X_k \prod_{c\neq b} \partial_c X_{\tau_c} &= \sum_{\substack{\tau \in S^n\\b=1\dots n}} \epsilon_\tau \partial_{\tau_b} v_{\tau_b} \partial_b X_{\tau_b} \prod_{c\neq b} \partial_c X_{\tau_c} + \\
    &\quad + \sum_{\substack{\tau \in S^n\\b=1\dots n}} \sum_{k \neq \tau_b} \epsilon_\tau \partial_k v_{\tau_b} \partial_b X_k \partial_{\tau^{-1}_k} X_k \prod_{c\neq b, \tau^{-1}_k} \partial_c X_{\tau_c}\\
    &= (\nabla \cdot v) \det( \nabla X ) +\\
    &\quad + \sum_{\substack{b=1\dots n\\k\neq b}} \sum_{\tau \in S^n} \epsilon_\tau \partial_{\tau_k} v_{\tau_b} \partial_b X_{\tau_k} \partial_k X_{\tau_k} \prod_{c\neq b, k} \partial_c X_{\tau_c}\\
    &=  (\nabla \cdot v) \det( \nabla X ) + 0
\end{align*}
The second term above is zero because replacing replacing $\tau$ with $\tau \circ (b\;\;k)$ in the inner sum produces a negative sign.

Similarly we have
$$
\sum_{\substack{\tau \in S^n\\b=1\dots n}} \epsilon_\tau \partial_k \sigma_{\tau_b,p} \partial_b X_k \prod_{c\neq b} \partial_c X_{\tau_c} \, dW_p = \partial_k \sigma_{k,p} \det( \nabla X )\, dW_p.
$$
For the last term in (\ref{D1}), the only difference is that we have a few extra cases to consider: When $l=\tau(b)$, $m = \tau(c)$, we will get $\det( \nabla X ) \partial_i \sigma_{ip} \partial_j \sigma_{jp}$. When $l = \tau(c)$ and $m = \tau(b)$, we will get $-\det( \nabla X) \partial_j \sigma_{ip} \partial_i \sigma_{jp}$. In all other cases we get $0$. This concludes proof of Lemma \ref{detAlemma}.

\begin{lemma}
For any time $t$, the map $X_t$ has a (spatially) smooth inverse.
\end{lemma}

\noindent{\bf Proof.} Define $\lambda$ by
$$
\lambda = \exp\left[ \int_0^t \left( \nabla \cdot v + 2 \nu E - \nu (\partial_k \sigma_{kp})^2 \right) \, dt  + \sqrt{2\nu} \int_0^t \partial_k \sigma_{kp} \,dW\sub{p}_s \right]
$$
The It\^o's formula immediately shows that $\lambda$ satisfies equation \eqref{deteq}. Since \eqref{deteq} is a linear SDE with smooth coefficients, uniqueness of the solution guarantees $D = \exp(\lambda)$ almost surely, and hence $D > 0$ almost surely.

The spatial invertibility of $X$ now follows as $X_t$ is locally orientation preserving and has degree $1$ (because $X_t$ is properly homotopic to $X_0$, the identity map). The (spatial) smoothness of the inverse is guaranteed by the inverse function theorem.\medskip

The above lemma shows existence of a spatial inverse of $X$. As before, we let $A$ denote the spatial inverse of $X$. We now derive a stochastic evolution equation of $A$ [equation \eqref{aeq}].

\begin{lemma}\label{l:itoComposed}
Let $Y$ be a $C^1$ stochastic flow of semi-martingales adapted to $\mathcal{F}_t$, the filtration of $W_t$. If for all $a \in \Rm^n$, $t > 0$ we have
$$
\int_0^t Y( X_s(a), ds ) = \int_0^t b( X_s(a), s)\,ds + \int_0^t \sigma'( X_s(a), s) dW_s
$$
then
$$
Y_t(a) = Y_0(a) + \int_0^t b(a, t) \,dt + \int_0^t \sigma'(a,t) \, dW_t.
$$
\end{lemma}

\noindent{\bf Proof.} Let $Y'$ be the process defined by
$$
Y'_t(a) = Y_0(a) + \int_0^t b(a, t) \,dt + \int_0^t \sigma'(a,t) \, dW_t,
$$
and set $\delta = Y - Y'$. Since $\delta$ is adapted to $\mathcal{F}_t$, there exists a non-negative predictable function $a$ such that
$$
\int_0^t a( x, y, s ) \, ds = \qv{\delta(x)}{\delta(y)}_t.
$$
Now, by definition of the generalized It\^o integral we have $\int_0^t \delta( X_s, s ) \equiv 0$ almost surely, and hence $\int_0^t a( X_s, X_s, s) \,ds \equiv 0$ almost surely. Since $X$ is a flow of homeomorphisms (diffeomorphisms actually), we must have $\forall t$, $a(x,x,t) \equiv 0$ almost surely. Thus $\delta = Y - Y'$ is of bounded variation.

Since we have shown above that $\delta$ has bounded variation,
$$
\int_0^t \delta( X_s, ds ) = \int_0^t \partial_t \delta\at{X_s, s} \, ds
$$
and hence $\forall t$, $\partial_t \delta_t \equiv 0$. At time $0$, $\delta_0 \equiv 0$ by definition, and hence $\delta_t \equiv 0$ almost surely for all $t$, concluding the proof.

\begin{lemma}\label{l:mpart-a}
There exists a process $B$ of bounded variation such that
\begin{equation}
\label{e:mpart-a} A_t = B_t - \sqrt{2\nu}\int_0^t (\nabla A_s) \sigma \, dW_s
\end{equation}
\end{lemma}
\noindent{\bf Proof.} Applying the generalized It\^o formula to $A \circ X$ we have
\begin{align}
\nonumber	0 &= \int_{t'}^t A(X_s, ds) + \int_{t'}^t \nabla A\at{X_s,s} \,dX_s + \tfrac{1}{2}\int_{t'}^t \partial^2_{ij} A\at{X_s, s} d\qv{X\sub{i}}{X\sub{j}}_s + \\
\nonumber	&\qquad\qquad + \qv{\int_{t'}^t \partial_i A(X_s, ds)}{X\sub{i}_t - X\sub{i}_{t'}}
\\
\label{e:daox}	&= \int_{t'}^t A(X_s, ds) + \int_{t'}^t \left[\nabla A\at{X_s,s} v + \nu a_{ij} \partial^2_{ij} A \at{X_s, s} \right]ds +\\
\nonumber	&\qquad\qquad + \sqrt{2\nu} \int_{t'}^t \nabla A\at{X_s,s} \sigma \,dW_s + \qv{\int_{t'}^t \partial_i A(X_s, ds)}{X\sub{i}_t - X\sub{i}_{t'}}.
\end{align}
Notice that the second and fourth terms on the right are of bounded variation. Applying Lemma \ref{l:itoComposed} we conclude the proof.

\begin{lemma}
The process $A$ satisfies the equation
\begin{equation}
\label{e:inverse-ito}	dA_t + (v \cdot \nabla) A_t \, dt - \nu a_{ij} \partial^2_{ij} A_t \, dt - \sqrt{2\nu} \partial_j A_t (\partial_i \sigma_{jk}) \sigma_{ik} \, dt + (\nabla A_t) \sigma \, dW_t = 0
\end{equation}
\end{lemma}

\noindent{\bf Proof.} Since the joint quadratic variation term in \eqref{e:daox} depends only on the martingale part of $\partial_i A$, we can compute it explicitly by
\begin{align}
\nonumber	\qv{ \int_{t'}^t \partial_i A( X_s, ds)}{X\sub{i}} &= -2\nu\int_{t'}^t \left( \partial^2_{ij} A_s \sigma_{jk} \sigma_{ik} + \partial_j A_s (\partial_i \sigma_{jk}) \sigma_{ik} \right)\circ X_s \,ds \\
\label{e:qv-A-X}
    &= -2\nu\int_{t'}^t \left( a_{ij} \partial^2_{ij} A_s + \partial_j A_s (\partial_i \sigma_{jk}) \sigma_{ik} \right)\circ X_s \,ds.
\end{align}

Substituting \eqref{e:qv-A-X} in \eqref{e:daox} and applying Lemma \ref{l:itoComposed} we conclude the proof.

\noindent{\bf{Acknowledgment.}}{ PC partially supported by NSF grant DMS-0504213. }

\end{document}